\documentclass[11pt]{article}
\usepackage[pdfpagelabels=true,plainpages=false,colorlinks=true,
linkcolor=blue,citecolor=blue,urlcolor=blue]{hyperref}
\usepackage[latin1]{inputenc}
\usepackage{algorithm}
\usepackage{algorithmic}
\usepackage{amsfonts}
\usepackage{amsmath}
\usepackage{amsthm}
\usepackage{cite}
\newcommand{\be}{\begin{equation}}
	\newcommand{\ee}{\end{equation}}
\newcommand{\bea}{\begin{eqnarray}}
	\newcommand{\eea}{\end{eqnarray}}
\usepackage{caption}
\usepackage{subcaption}

\usepackage{graphicx,epsfig}
\usepackage{enumerate}
\usepackage{setspace}
\usepackage{amssymb}
\usepackage{amsmath}
\usepackage{amsthm}
\usepackage{appendix}
\newtheorem{alg1}{\textbf{Algorithm}}[section]
\numberwithin{equation}{section}
\newtheorem{thm}{Theorem}[section]
\theoremstyle{definition}
\newtheorem{dfn}{Definition}[section]
\theoremstyle{Remark}
\theoremstyle{lemma}

\newtheorem{lem}{Lemma}[section]

\theoremstyle{example}
\begin{document}

\title{Global Convergence and Rate Analysis of the Steepest Descent Method for Uncertain Multiobjective Optimization via a Robust Optimization Approach}
\author{Shubham Kumar\footnote{Roorkee Institute of Technology, Roorkee, India, E-mail: kumarshubham3394@gmail.com, corresponding author},
	Nihar Kumar Mahato \footnote{IIITDM Jabalpur, India, E-mail: nihar@iiitdmj.ac.in},
	Debdas Ghosh \footnote{IIT(BHU), India, E-mail: debdas.mat@iitbhu.ac.in}
	\\}
\date{}
\maketitle
\begin{abstract}
In this article, we extend our previous work (Applicable Analysis, 2024, pp. 1-25) on the steepest descent method for uncertain multiobjective optimization problems. While that study established local convergence, it did not address global convergence and the rate of convergence of the steepest descent algorithm. To bridge this gap, we provide rigorous proofs for both global convergence and the linear convergence rate of the steepest descent algorithm. Global convergence analysis strengthens the theoretical foundation of the steepest descent method for uncertain multiobjective optimization problems, offering deeper insights into its efficiency and robustness across a broader class of optimization problems. These findings enhance the method's practical applicability and contribute to the advancement of robust optimization techniques.
\end{abstract}

\section{Introduction}
\label{intro}
Multiobjective optimization problems (MOP) are essential in decision-making processes where multiple conflicting objectives need to be optimized simultaneously. Unlike single-objective optimization problems, which focus on identifying a single best solution, MOP aims to determine a set of solutions known as the Pareto optimal set. In this set, enhancing one objective is only possible at the expense of deteriorating another. MOP is vital in real-world optimization applications across various domains, including business, economics, and diverse scientific and engineering disciplines (see, for example, \cite{bhaskar2000applications, stewart2008real, miettinen1999nonlinear, deb2013multi} and the references cited therein). Various researchers have contributed to the development of MOP results in various contexts (e.g., \cite{ehrgott2005multicriteria, miettinen1999nonlinear, fliege2000steepest, fliege2009newton, fukuda2014survey, ansary2015modified, lai2020q} and the references cited therein). Many optimization methods in the literature commonly assume precise values for the parameters involved. However, real-world scenarios often contain inherent uncertainties in data or model parameters due to estimation errors, rounding discrepancies, or implementation challenges. To address such uncertainties, various probabilistic and non-probabilistic approaches have been developed, including stochastic programming \cite{birge2011introduction} and robust optimization \cite{doolittle2018robust, chen2012including, yu2013robust, fliege2014robust}. In stochastic programming, uncertain coefficients are modeled as random variables with known distribution functions, which can be difficult to determine in optimization problems. Interval optimization handles uncertainty by representing uncertain parameters as intervals. In contrast, robust optimization does not rely on probability distributions, making it more resilient to errors in estimating uncertainties. Additionally, it does not require uncertainty to be expressed in interval form. 
Robust optimization seeks to identify solutions that remain effective across a range of possible future scenarios. Typically, it transforms the original problem into a deterministic equivalent, referred to as the robust counterpart, enabling the application of standard optimization algorithms to obtain solutions that are resilient to uncertainties. Such problems commonly arise in various contemporary research fields (see, for example, \cite{dellnitz2009computation, witting2013variational, gunawan2005multi, barrico2006robustness, goerigk2015algorithm, doolittle2018robust} and the references cited therein).
Jeyakumar et al. \cite{jeyakumar2012robust} introduced a robust duality theory for generalized convex programming problems involving data uncertainty within the robust optimization framework. Goberna et al. \cite{goberna2013robust} formulated a duality theory for semi-infinite linear programming problems under uncertainty in constraint functions, the objective function, or both. Ide et al. \cite{ide2014relationship} investigated the relationship between different robustness concepts for MOP, established connections with set order relations, and developed algorithms for solving uncertain multiobjective optimization problems (UMOP). Ehrgott et al. \cite{ ehrgott2014minmax} studied a class of UMOPs and proposed methods to determine both Pareto and weak Pareto optimal solutions. Goberna et al. \cite{goberna2015robust} analyzed robust solutions for multiobjective linear programs with uncertainties in both objective functions and constraints. Chuong \cite{ chuong2016optimality} derived optimality and duality results for robust multiobjective optimization problems, while Lee and Lee \cite{lee2018optimality} established optimality conditions and duality theorems for robust semi-infinite multiobjective optimization problems.
\par Recently, in the field of numerical optimization, we developed the steepest descent method \cite{kumar2024steepest}, Newton's method \cite{shubham2023newton}, for uncertain multiobjective optimization problems under finite uncertainty sets. Later, we proposed a modified quasi-Newton method \cite{kumar2024modified} to address uncertain multiobjective optimization problems under a finite uncertainty set.
\par In \cite{kumar2024steepest}, we have extended the steepest descent method from deterministic multiobjective optimization problems to UMOP. In \cite{kumar2024steepest}, a steepest descent algorithm was introduced for solving the objective-wise worst-case-type robust counterpart (OWRC) of the UMOP, along with a discussion on the convergence analysis of the sequence generated by the proposed algorithm. Although a numerical experiment was conducted, the authors did not demonstrate the full convergence of the steepest descent method. Additionally, there was no discussion on the convergence rate of the steepest descent method for OWRC. Motivated by this, the present study aims to establish the full convergence of the steepest descent method for OWRC in addressing the solution of UMOP. Furthermore, we prove that the sequence generated by the steepest descent algorithm converges to the solution with a linear rate. 
\par To the best of our knowledge, no literature discusses the complete convergence of the steepest descent algorithm for OWRC in solving UMOP and presents its linear convergence rate. Therefore, the main objectives of the present study are to develop the full convergence analysis of the steepest descent algorithm and to establish its linear convergence rate.
\par The structure of this article is as follows: In Section \ref{sec2}, we introduce important notations and results relevant to our problems. Section \ref{sec3} provides a brief overview of the steepest descent method for OWRC to find the solution of $P(U).$ In Section \ref{sec4}, we establish the full convergence of the steepest descent algorithm. Section \ref{sec5} focuses on establishing the linear convergence rate of the steepest descent algorithm. Finally, Section \ref{sec6} concludes with some remarks. 

\section{Preliminaries}\label{sec2}
Throughout the article, we use the following notations: The set of real numbers is denoted by $\mathbb{R}$. The sets of non-negative and positive real numbers are represented by $\mathbb{R}_{\geq}$ and $\mathbb{R}_{>}$, respectively. The $n$-dimensional space of real-valued tuples is denoted by $\mathbb{R}^n$, while $\mathbb{R}_{\geq}^n$ and $\mathbb{R}_{>}^n$ represent spaces where each component is non-negative and positive, respectively. The notation $\bar \Lambda$ denotes the set $\{1, 2, \dots, p\}$, and $\Lambda$  represents the set $\{1, 2, \dots, m\}$. For vectors $s, q \in \mathbb{R}^n$, the relation $s \leq q$ means each component of $s$ is less than or equal to the corresponding component of $q$. Similarly, $s\geq q$ means each component of $s$ is greater than or equal to the corresponding element of $q$. The notation $s > q$ signifies that each component of $s$ is strictly greater than the corresponding component of $q$.

Two key challenges often hinder the application of optimization techniques in real-world scenarios. First, handling input data that is not precisely known poses a significant issue in almost all practical optimization problems. This uncertainty may arise due to measurement errors, data inaccuracies, future uncertainties, fluctuations, or external disturbances. Second, real-world optimization problems typically involve multiple decision-makers and diverse objectives. As a result, they rarely have a single well-defined objective function but instead require consideration of multiple competing goals. Before starting uncertain multiobjective optimization, we start introducing the multiobjective optimization problems. 
\par An unconstrained deterministic multiobjective optimization problem (DMOP) is given by  
\begin{equation*}
\rm{DMOP}:~\min_{\textit{x}\in\mathbb{R}^n}F(\textit{x}),
\end{equation*}
where \( F:\mathbb{R}^n\to\mathbb{R}^m \) is defined as  
\begin{equation*}
F(x) = (F_1(x), F_2(x), \dots, F_m(x)),
\end{equation*}
with each \( F_j:\mathbb{R}^n\to\mathbb{R} \) for \( j = 1,2,\dots,m \). The solution of a DMOP can be defined in terms of a Pareto optimal solution (POS) 
\big(equivalently, an efficient solution (ES)\big) and a weak Pareto optimal solution (WPOS) 
\big(also known as a weak efficient solution (WES)\big).

\begin{dfn} Given a point $x^* \in D\subset \mathbb{R}^n$,
	\begin{enumerate}
		\item[(i)] $x^*$ is said to be an ES or POS of $F$ if and only if there does not exist $x \in D$ such that $F(x) \leq F(x^*)$ and $F(x) \neq F(x^*)$;
		\item[(ii)] $x^*$ is said to be WES or WPOS of $F$ if and only if there does not exist $a \in D$ such that $F(x) < F(x^*)$;
		\item[(iii)] $x^*$ is said to be locally ES or locally WES $F$ if and only if there exists a neighborhood $X \subset D$ such that the point $x^*$ is ES or WES for $F$ restricted to $X$.
	\end{enumerate}
\end{dfn}
If $x^*$ is an efficient solution, then $F(x^*)$ is called a non-dominated point, and the set of efficient solutions and the set of non-dominated points (Pareto front) are called the efficient set and the non-dominated set, respectively. For further details of these definitions, see references 
\cite{ehrgott2005multicriteria, miettinen1999nonlinear,fliege2000steepest,drummond2005steepest,fukuda2014survey}. Some scalarization methods to solve POS have been developed. For an overview of scalarization methods, refer to Ehrgott \cite{ehrgott2005multicriteria}, and Miettinen \cite{miettinen1999nonlinear}. Apart from scalarization methods, several classical derivative-based methods for scalar optimization have been extended to multiobjective optimization. These derivative-based methods for multiobjective optimization are as follows: steepest descent method \cite{fliege2000steepest,drummond2005steepest}; Newton's method \cite{fliege2009newton}; quasi-Newton method \cite{lai2020q,ansary2015modified,povalej2014quasi,mahdavi2020superlinearly,morovati2018quasi,qu2011quasi}; Conjugate gradient method\cite{lucambio2018nonlinear,gonccalves2020extension}; Projected gradient method \cite{drummond2004projected,fukuda2011convergence,fukuda2013inexact,fazzio2019convergence}; projected gradient method [1]; and proximal gradient method \cite{bonnel2005proximal,ceng2010hybrid}. Apart from DMOP, we define an uncertain multiobjective optimization problem as follows:
\par An uncertain multiobjective optimization problem with finite uncertainty set $U,$ is defined as:  
\begin{equation}\label{mp}
	P(U)=\{P(\xi_i):i\in\bar{\Lambda}\}, ~U=\{\xi_i:i\in\bar{\Lambda}\},
\end{equation}
where for any fix $\xi_i\in U\subset \mathbb{R}^k$, $$P(\xi_i):\displaystyle \min_{x\in \mathbb{R}^n} f(x,\xi_i)$$ and $f:\mathbb{R}^n\times U\to \mathbb{R}^m$ such that $f(x,\xi_i)=( f_1(x,\xi_i),f_2(x,\xi_i),...,f_m(x,\xi_i)), ~~ \xi_i\in U\subset \mathbb{R}^k.$\\
\par The main objective for the problem $P(U)$ is to find a robust Pareto optimal solution (RPOS) or a robust weakly Pareto optimal solution (RWPOS), which can be defined as follows:
\begin{dfn}\cite{ehrgott2014minmax} Given a point $x^* \in D\subset \mathbb{R}^n$,
	\begin{enumerate}
		\item[(i)] $x^*\in D$ is said to be RPOS for $P(U)$ if and only if there is no $x\in D-\{x^*\}$ such that $f(x;U)\subset f(x^*;U)-\mathbb{R}^k_\ge$;
		\item[(ii)] $x^*\in D$ is said to be RWPOS for $P(U)$ if and only if there is no  $x\in D-\{x^*\}$ such that $f(x;U)\subset f(x^*;U)-\mathbb{R}^k_>,$
		where at any $x\in D,$ $f(x; U)=\{f(x,\xi_i):\xi_i\in U\}$ i.e., $f(x; U)$ is the set of all possible values of the objective function at $x.$ 
	\end{enumerate}
\end{dfn}
Now, we define the objective-wise worst-case cost-type robust counterpart (OWRC) of 
$P(U)$ as follows:
\begin{equation}\label{1.2}
		\rm{OWRC:}~~\min_{\textit{x}\in \mathbb{R}^n}\Phi(\textit{x}),
\end{equation} where $\Phi(x)=(\Phi_1(x),\Phi_2(x),\ldots,\Phi_m(x)),$ and $\Phi_j(x)=\displaystyle\max_{i\in\bar{\Lambda}}f_j(x,\xi_i).$\\
Ehrgott et al. \cite{ehrgott2014minmax} demonstrated that the solution to OWRC will also solve $P(U).$ Therefore, solving OWRC offers an advantage over solving P(U), as P(U) represents a collection of problems whereas OWRC is a single deterministic multiobjective problem. 
Now, we present a modified version of Theorem 4.11 (given in \cite{ehrgott2014minmax}) for $P(U)$ and OWRC, establishing the connection between the solutions of OWRC and $P(U)$.
\begin{thm}\label{t1} \cite{kumar2024steepest} Let $P(U)$ be an UMOP with a finite uncertainty set, and OWRC be the robust counterpart of $P(U)$. Then:
	\begin{enumerate}[(a)]
		\item If $x^*\in\mathbb{R}^n$ is a POS to $OWRC$, then $x^*$ is RPOS for $P(U)$.
		\item If $\displaystyle\max_{\xi\in U} f_j(x,\xi)$  exist for all $j\in\Lambda$ and all $x \in \mathbb{R}^n$ and $x^*$ is a WPOS $OWRC,$ then $x^*$ is a robust weakly Pareto optimal solution for $P(U)$.
	\end{enumerate}
\end{thm}
\par Note that \( \Phi_j(x) = \displaystyle\max_{i \in \bar{\Lambda}}f_j(x, \xi_i) \) is generally not a differentiable function, thereby making OWRC a deterministic non-smooth MOP. Several versions of the steepest descent method for smooth MOPs have been developed in the literature; for further details, see \cite{fliege2000steepest,drummond2005steepest,fukuda2014survey}. Additionally, let's define $I_j(x)$ as the set of all points in $\bar{\Lambda}$ where $f_j(x, \xi_i)$ achieves its maximum value. In other words, $I_j(x)$ is given by $I_j(x) = \{i \in \bar{\Lambda} : f_j(x, \xi_i) = \Phi_j(x)\}$. Further details can be found in \cite{sun2006optimization}. 
\par Based on Theorem \ref{t1}, the critical point of OWRC serves as the robust critical point of $P(U)$, which is defined as follows:
\begin{dfn}{(\textbf{Critical point for $OWRC$ or Robust Pareto critical point for $P(U)$})}\cite{kumar2024steepest}\\
	Let  $ \Phi:\mathbb{R}^n \to \mathbb{R}^m$ be defined as $\Phi(x)=(\Phi_1(x),\Phi_2(x),\ldots,\Phi_m(x))$, where $\Phi_j(x)= \displaystyle \max_{i\in \bar{\Lambda}} f_j(x,\xi_i),$ $j\in \Lambda$ and $f_j:\mathbb{R}^n \times U \to \mathbb{R}$ is a continuously differentiable for each $j\in\Lambda$ and $\xi_i\in U$. Then, a point $x^*$ is called a critical point for $\Phi$ (for OWWC) if for all $v\in\mathbb{R}^n$, there exists a  $j^0\in \Lambda$ such that $\nabla f_{j^0}(x^*,\xi_i)^Tv\geq0$. In other words, if  $x^*$ is  a critical point for problem $OWWC$, then there does not exist $v\in \mathbb{R}^n$  such that $ \nabla f_j(x,\xi_i)^Tv <0,~ \forall ~i \in I_j(x)$ and $j \in \Lambda $.
	Also, if  $x^*$ is a critical point for problem OWRC then it will be a robust Pareto critical point for $P(U).$
\end{dfn} 
\begin{dfn}
	Let $ \digamma:\mathbb{R}^n \to \mathbb{R}$ be a function. Then a vector $v\in \mathbb{R}^n$ is said to be descent direction of $h$ at $x$ if and only if there exists $\delta > 0$ such that $\digamma(x + \alpha v) < \digamma(x),$ $\forall\alpha\in (0, \delta).$ If $h$ is continuously differentiable then a vector $v\in\mathbb{R}^n$ is said to be a descent direction for $h$ at $x$ if $\nabla \digamma(a)^Tv< 0$.\\
	In particular, for OWRC a vector $d\in\mathbb{R}^n$ is said to be descent direction for $T$ at $x$ if
	\begin{center}
		$	 \nabla f_j(x,\xi_i)^Tv< 0, ~\forall j \in \Lambda$ and $i\in I_j(x).$
	\end{center}
	Also, if $v\in\mathbb{R}^n$ is a descent direction for $\Phi$ at $x$, then there exists $\bar \alpha >0$ such that \begin{center} $\Phi_j(x+\alpha v)< \Phi_j(x), ~\forall~ j \in \Lambda$, ~and $\forall \alpha \in(0,\bar{\alpha}]$.\end{center}
\end{dfn}
\begin{thm}\label{Thm 2.1}\cite{kumar2024steepest}
	Let  $f_j(x,\xi_i)$ is continuously differentiable  for each $j\in \Lambda$ and $i \in \bar{\Lambda}$. If $x^*$ is a weakly efficient solution for OWRC, then $0 \in conv \left( \displaystyle \cup_{j=1}^{m}\partial \Phi_j(x^*) \right).$
\end{thm}
\begin{thm}\label{Thm 2.301}\cite{kumar2024steepest}
	If  $f_j(x,\xi)$ is continuously differentiable and convex for each $j\in \Lambda$ and every $\xi\in U$, then  $x^*$ is a weakly efficient solution for OWRC if and only if
	\[0 \in conv \left( \displaystyle \cup_{j=1}^{m}\partial \Phi_j(x^*) \right).\]
\end{thm}
\begin{dfn}\label{qfd} \cite{fliege2000steepest}
	A sequence $\{b^k\}\subset \mathbb{R}^n$ is quasi-Fejer convergent to a set $B\subset\mathbb{R}^n$ if for every $b\in B$ there exists a sequence $\{\delta_k\}\subset \mathbb{R},$ $\delta_k\geq0$ for all $k$ and such that $\|b^{k+1}-b\|^2\leq\|b^k-b\|^2+\delta_k$ for all $k=0,1,\ldots$ with $\sum_{k=0}^{\infty} \delta_k<\infty.$
\end{dfn}
\begin{thm}\label{qfc}\cite{fliege2000steepest}
	If a sequence $\{b^k\}$ is quasi-Fejer convergent to a non empty set $B\subset\mathbb{R}^n,$ then $\{b^k\}$ is bounded. Furthermore, if $\{b^k\}$ has a cluster point $b$ which belongs to $B,$ then $\lim_{k\to \infty}b^k=b.$ 
\end{thm}

\section{The steepest descent method for OWRC}\label{sec3}
The steepest descent method helps to find the critical point for OWRC, which will be the robust Pareto critical point for $P(U).$~Let $ \Phi:\mathbb{R}^n \to \mathbb{R}^m$ be defined as $\Phi(a)=(\Phi_1(x),\Phi_2(x),\ldots,\Phi_m(x))$, where $\Phi_j(x)= \displaystyle \max_{i\in \bar{\Lambda}} f_j(x,\xi_i),~ j\in \Lambda$ and $f_j:\mathbb{R}^n \times U \to \mathbb{R}$ is a continuously differentiable and convex function for each $x$ and $\xi_i\in U$. 
Two tools are required to develop the steepest descent algorithm: one is the steepest descent direction, and the second is the step length size in the steepest descent direction. 
\subsection{A subproblem to find descent directions for OWRC}\label{ss3.1}
In \cite{shubham2024steepest}, to find the steepest descent direction for OWRC, we have constructed the following unconstrained real-valued minimization subproblem.
\begin{equation}\label{eq31}
	\displaystyle \min_{t\in \mathbb{R}^n} \vartheta_{x}(t)+\frac{1}{2}\|t\|^2,
	\text{ where }\vartheta_{x}(t)= \displaystyle \max_{j\in\Lambda}\max_{i\in \bar{\Lambda}}  \{ f_j(x,\xi_i)+\nabla f_j(x,\xi_i)^Tt-\Phi_j(x)\}.
\end{equation}
 Note that being the maximum of a max-linear function, as $\vartheta_{x}(t)$ is convex, therefore the given objective function in the subproblem \eqref{eq31} is strongly convex. Hence the subproblem \eqref{eq31} will have a unique solution, which will be the steepest descent direction for OWRC. Let $t(x)$ and $\Theta(x)$ be the optimal solution and optimal value of the subproblem \eqref{eq31}, respectively. Then,
\begin{eqnarray}
	t(x)&=&\underset{v\in R^n}{\mathrm{argmin}}  \left( \vartheta_{x}(t)+\tfrac{1}{2}\|t\|^2 \right)=	-	\sum\limits_{j\in \Lambda}  \sum\limits_{i\in \bar \Lambda } \lambda_{ij} \nabla f_j(x,\xi_i),\label{ov}\\
	\text{and} ~~~~~\Theta(x)&=&  \vartheta_{x}(t(x))+\frac{1}{2}\|t(x)\|^2\label{os}.
\end{eqnarray}
Thus, $t(x)$ is the solution of $P(x)$ and ultimately, it will also be the solution of the subproblem \eqref{eq31}.
\par In the following theorem, the characterization of a critical point in terms of $v(x)$ and $\Theta(x)$ is presented. The properties of $t(x)$ and $\Theta(x)$ are discussed, along with the proof that $t(x)$ serves as a descent direction for $\Phi$ at $x$
\begin{thm}\label{Thm 3.2.} Let $t(x)$  and $\Theta(x)$ be the solution  and optimal value of the subproblem \eqref{eq31}, respectively, which are given in the equations $(\ref{ov})$ and $(\ref{os})$. Then, the following results will hold:
	\begin{enumerate}
		\item $ t(x)$ is bounded on any compact subset $C$ of  $\mathbb{R}^n.$
		\item The following conditions are equivalent:
		\begin{enumerate}
			\item The point $x$ is a noncritical point.
			\item $\Theta(x)<0.$
			\item $t(x)\not =0.$
			\item$t(x)$ is a descent direction for $\Phi$ at $x$.
		\end{enumerate}
	\end{enumerate}
	In particular, $x$ is a  critical point $\iff$ $\Theta(x)=0$.
\end{thm}
\begin{proof}
	See the Theorem 3.2 in \cite{shubham2024steepest}.
\end{proof}
To find the step length size following Armijo-type inexact line search technique is used in \cite{shubham2024steepest}. 
\begin{thm}\label{3.2} Let $x$ be any non-critical point. Then, for any $\beta > 0$ and $\epsilon\in(0,1]$ there exists
	an $\alpha \in [0,\epsilon]$ such that
	\begin{equation}\label{13}
		\Phi_j(x+\alpha t)\leq  \Phi_j(x)+  \alpha \beta \Phi^*_j(x,t),~
		\text{where}~\Phi^*_j(x,t)= \max_{i\in\bar\Lambda}\{ f_j(x,\xi_i) + \nabla f_j(x,\xi_i)^Tv \}-\Phi_j(x),~ j\in \Lambda.
	\end{equation}
\end{thm}
\begin{proof}
	For the complete proof, see the subsection 3.2 in \cite{shubham2024steepest}. 
\end{proof}
Equation (\ref{13}) represents the step size rule for the steepest descent algorithm for the OWWC.
\par Utilizing the steepest descent direction and an inexact line search technique to determine the step length, the steepest descent algorithm for OWRC is as follows.
\begin{alg1}\label{algo1}(\textbf{Steepest descent algorithm for OWRC})
	\begin{enumerate}[{Step} 1]
			\item Choose $\epsilon>0$, $\beta \in (0,1) $   and $x^0\in \mathbb{R}^n$. Set $k :=0.$
		\item Solve subproblem \ref{eq31} and find $t^k$, $\Theta(x^k)$.
		\item If  $|\Theta (x^k)|<\epsilon$, then stop. Otherwise, proceed to Step \ref{st4}.
		\item \label{st4} Choose $\alpha_k$ as the largest $\alpha \in \{ \frac{1}{2^r} : r=1,2,3,\ldots\} $ satisfying (\ref{13}).
		\item Define $x^{k+1}:= x^k + \alpha_{k} t^k$, update $k:=k+1$ and go to Step 2.
	\end{enumerate}
\end{alg1}

\begin{thm}\label{3.4}
	Let $\{x^k\}$ be a sequence that is produced by Algorithm \ref{algo1}. Assume that $ t(x)$ is bounded on every compact subset $C$ of ~ $\mathbb{R}^n\times \mathbb{R}^k$. Then, any accumulation point of $\{x^k\}$ is a critical point for $\Phi$.
\end{thm}
\begin{proof}
	See the Theorem 3.3 in \cite{shubham2024steepest}.
\end{proof}
In the upcoming section, we establish the global convergence of the steepest descent method.
\section{Global convergence of steepest descent method}\label{sec4}
\begin{lem}\label{gl1}
	Let $f_j:\mathbb{R}^n\times U\to \mathbb{R}$ be continuously differentiable and $d\in\mathbb{R}^n$ such that $\nabla f_j(a,\xi_i)^Tt<0$ for each $i$ and $j$. Then for any $\beta\in(0,1)$ there exists some $\epsilon>0$ (depending on $x,~t,$ and $\beta$) such that 
	\begin{equation}
		f_j(a+\alpha t,\xi_i)< f_j(a,\xi_i)+\beta\alpha \nabla f_j(a,\xi_i)^T t,~~~~~\forall\alpha\in(0,\epsilon].
	\end{equation} 
\end{lem}
\begin{proof}
	Since $f_j(a,\xi_i)$ is continuously differentiable and $\nabla f_j(a,\xi_i)^Tt<0$ for each $i$ and $j.$ Then, we obtained 
	\begin{equation*}
		\displaystyle\lim_{\alpha\to \infty}\frac{ f_j(a+\alpha t, \xi_i)-f_j(a,\xi_i)}{t}=\nabla f_j(a,\xi_i)^Tt.
	\end{equation*}
	As $\beta\in(0,1)$ and $\nabla f_j(a,\xi_i)^Tt<0$ for each $i$ and $j,$ therefore there must exists $\epsilon>0$ such that 
	\begin{equation}
		f_j(a+\alpha d,\xi_i)- f_j(a,\xi_i)<\beta\alpha \nabla f_j(a,\xi_i)^T t,~~~~~\forall\alpha\in(0,\epsilon],
	\end{equation}
	which is the required inequality.
\end{proof}
\begin{lem}\label{gl2} Let $f_j:\mathbb{R}^n\times U\to \mathbb{R}$ be continuously differentiable
	and $\tilde{x}$ be a weak Pareto optimal solution for $\Phi.$ Then $\tilde{x}$ is a critical point. Moreover, if $f_j(x,\xi_i)$ is convex, then the critical point will be a weak Pareto optimal point. Therefore, Pareto optimality and criticality are equivalent under the convexity of objective functions.
\end{lem}
\begin{proof} It is given that $\tilde{x}$ is a a WPOS for $\Phi.$ On contrary, we assume $\tilde{x}$ is not a critical point, then there must exist a $d\in \mathbb{R}^n$ such that \begin{equation}\label{gl2e1}
		\nabla f_j(x,\xi_i)^Tt\leq0,~ \text{for each} ~i ~\text{and}~ j.
	\end{equation} Now, by Lemma \ref{gl1}, we obtain 
	\begin{equation}\label{gl2e2}
		f_j(\tilde{x}+\alpha t,\xi_i)\leq f_j(\tilde{x},\xi_i)+\alpha \nabla f_j(\tilde{x},\xi_i)^Tt ~\text{for each} ~i ~\text{and}~ j.
	\end{equation} Therefore, by inequalities (\ref{gl2e1}) and (\ref{gl2e1}), we obtain
	$$f_j(\tilde{x}+\alpha t,\xi_i)\leq f_j(\tilde{x},\xi_i),~\text{for each} ~i ~\text{and}~ j,$$
	which implies $$\max_{i\in\bar{\Lambda}}f_j(\tilde{x}+\alpha t,\xi_i)\leq\max_{i\in\bar{\Lambda}} f_j(\tilde{x},\xi_i),$$ therefore $$\Phi(\tilde{x}+\alpha t)\leq \Phi(\tilde{x}).$$ The last inequality contradicts the fact that $\tilde{x}$ is a WPOS for $\Phi.$ Hence $\tilde{x}$ is a critical point.\\
	Now if we consider $f_j(x,\xi_i)$ is convex for each $i$ and $j.$ It is assumed that $\tilde{x}$ is a critical point. We have to show that $\tilde{x}$ is a WPOS. By definition of critical point there exists $\tilde{j}= j(\tilde{x}) \in\Lambda$ such that \begin{equation}\label{gl2e3}
		\nabla f_{\tilde{j}}(\tilde{x},\xi_i)^Tt\geq0,~\forall d\in \mathbb{R}^n,~i\in\bar{\Lambda}.
	\end{equation} By the convexity of $f_{\tilde{j}}(x,\xi_i),$ we have 
	\begin{equation}\label{gl2e4}
		f_{\tilde{j}}(y,\xi_i)-f_{\tilde{j}}(x,\xi_i)\geq \nabla f_{\tilde{j}}(x,\xi_i)^T(y-x)~ \forall x,y\in \mathbb{R}^n,~i\in \bar{\Lambda},~j\in\Lambda.
	\end{equation} 
	Now if we put $x=\tilde{x}$ and take $d=y-x$ then by (\ref{gl2e3}) and (\ref{gl2e4}), we obtained 	$$f_j(y,\xi_i)\geq f_j(\tilde{x},\xi_i),~\text{for each}~ i\in \bar{\Lambda}~\text{and}~\tilde{j}\in\Lambda,$$ which implies $$\max_{i\in\bar{\Lambda}}f_j(y,\xi_i)\geq \max_{i\in\bar{\Lambda}}f_j(\tilde{x},\xi_i),~\tilde{j}\in\Lambda,$$
	therefore $$\Phi_{\tilde{j}}(y)\geq \Phi_{\tilde{j}}(\tilde{x}),~\tilde{j}\in\Lambda,$$
	and hence $\tilde{x}$ is a weak Pareto optimal point.
	
\end{proof}
\textbf{Assumption 1:} Every sequence $\{z^k\}\subset \Phi(\mathbb{R}^n)=\{f(x):x\in \mathbb{R}^n\}$ which is decreasing in component-wise sense is bounded below by a point, $\Phi(\mathbb{R}^n)$ i.e., for any $\{z^k\}$ contained in $\Phi(\mathbb{R}^n)$ with $z^{k+1}<z^k,$ for all $k,$ there exists $\hat{z}$ such that $\Phi(\hat{z})\leq z^k$ for all $k.$\\
\begin{lem}\label{4.3}
	Suppose that $f_j(x,\xi_i)$ is continuously differentiable and convex for each $j,$ $i$ and let $\tilde{x}$ be such that $f_j(\tilde{x},\xi_i)\leq f_j(x^k,\xi_i)$ for some $k.$ Then, the following inequality holds
	\begin{equation}
		\|\tilde{x}-x^{k+1}\|^2\leq 	\|\tilde{x}-x^{k}\|^2+ 	\|x^k-x^{k+1}\|^2.
	\end{equation}
\end{lem}
\begin{proof}
	Since $f_j(x,\xi_i)$ is continuously differentiable and convex, we have \begin{equation*}
		f_{\tilde{j}}(\tilde{x},\xi_i)-f_{\tilde{j}}(x^k,\xi_i)\geq \nabla f_{\tilde{j}}(x^k,\xi_i)^T(\tilde{x}-x^k),~\text{for each}~i\in \bar{\Lambda},~j\in\Lambda. 
	\end{equation*} 
	It is also given that 
 $f_j(\tilde{x},\xi_i)\leq f_j(x^k,\xi_i)$ then by the above inequality we get \begin{equation}\label{3.8}
		\nabla f_{\tilde{j}}(x^k,\xi_i)^T(\tilde{x}-x^k)\leq0~\text{for each}~i,~j. 
	\end{equation}
	Now by equation (\ref{4***}), we have 	$t(x^k) = -	\sum\limits_{j\in \Lambda}  \sum\limits_{i\in \bar \Lambda } \lambda_{ij} \nabla f_j(x^k,\xi_i).$ Also, 
	\begin{align}\label{3.9}
		-t(x^k)(\tilde{x}-x^k)= 	\sum\limits_{j\in \Lambda}  \sum\limits_{i\in \bar \Lambda } \lambda_{ij} \nabla f_j(x^k,\xi_i)(\tilde{x}-x^k).
	\end{align}
	Then, by (\ref{3.8}) and (\ref{3.9}), we obtain
	\begin{align}\label{3.10}
		-t(x^k)(\tilde{x}-x^k)\leq0.
	\end{align}
	Therefore, from the above inequality (\ref{3.10}) we get 
	\begin{align}\label{3.11}
		(x^k-x^{k+1})(\tilde{x}-x^k)\leq0~~~~~~~~~~~~\because~ x^k-x^{k+1}=-\alpha_k t^k~\text{with}~\alpha_k>0.
	\end{align}
	Thus, $$\|\tilde{x}-x^{k+1}\|^2=\|\tilde{x}-x^k\|^2+\|x^k-x^{k+1}\|+2(x^k-x^{k+1})^T(\tilde{x}-x^k),$$
	hence by the above inequality and (\ref{3.11}), we get the required inequality.
	\begin{lem} \label{4.4} Let $\{x^k\}$ be any infinite sequence of non-critical point in $\mathbb{R}^n$ such that $\Phi(x^{k+1})\leq \Phi(x^k)$ for all $k.$
		\begin{enumerate}[(i)]
			\item \label{item1} If $x^*$ is any accumulation point of $\{x^k\}$, then $\Phi(x^*)\leq \Phi(x^k)$ for all $k$ and $\displaystyle\lim_{k\to \infty}\Phi(x^k)=\Phi(x^*).$ Moreover, $\Phi$ is constant in the set of accumulation points of $\{x^k\}.$
			\item If $\{x^k\}$ is generated by Algorithm \ref{algo1} and has an accumulation point, then all the results of item (\ref{item1}) hold for this sequence.
		\end{enumerate}
	\end{lem}
	\begin{proof}
		Let $x^*$ be any accumulation point of  $\{x^{k}\}.$ Then there must exists a subsequence $\{x^{k_l}\}_l$ of  $\{x^{k}\}$ that converges to $x^*.$ Now take a fixed but arbitrary $k,$ and for $l$ large enough, we have $k_l>k$ and \begin{equation}\Phi(x^{k_l})\leq \Phi(x^k).
		\end{equation} As $l \to \infty$ we obtain \begin{equation}\label{3.122}
			\Phi(x^*)\leq \Phi(x^k),~\text{for all}~k,\end{equation} 
		and in the component-wise sense we can write this inequality as \begin{equation}\label{3.133}
			\Phi_j(x^*)\leq \Phi_j(x^k).
		\end{equation}
		Also,  \begin{equation}\label{3.14}\Phi_j(x^{k_l})\to \Phi_j(x^*)~ \text{as}~ l\to \infty~ \text{for all}~j.\end{equation}
		Therefore, by (\ref{3.133}) and (\ref{3.14}), we obtain \begin{equation}\label{3.15}
			\Phi(x^k)\to \Phi(x^*)~ \text{as} ~k\to \infty.
		\end{equation}
		\\To prove that the $\Phi$ is constant in the set of accumulation points, we assume $\bar{x}$ is any other accumulation point of $\{x^k\}.$ Then there must exists a subsequence $\{x^{k_r}\}$ of $\{x^k\}$ such that $x^{k_r}\to\bar{x}$ as $r\to \infty.$ Now, by (\ref{3.122}), $\Phi(x^*)\leq \Phi(x^{k_r}),~\text{for all}~r,$ therefore as $r\to \infty$, we obtain $\Phi(x^*)\leq \Phi(\bar{x}).$ Also by interchanging the $x^*$ and $\bar{x}$ we get $\Phi(\bar{x})\leq \Phi(x^*).$ Then by last two inequality we get $\Phi(\bar{x})= \Phi(x^*)$ and hence $\Phi$ is constant. \\
		By Algorithm 1, we can observe that the sequence $\{x^k\}$ generated by Algorithm 1 satisfies $\Phi(x^{k+1})\leq \Phi(x^k)$ therefore all the results of item (ii) automatically hold.
	\end{proof}
	\begin{lem}\label{4.5}
		Suppose that $\{\Phi(x^k)\}$ is component wise bounded below by $\hat{y}\in \mathbb{R}^m.$ Then
		\begin{equation}
			\sum_{k=0}^{\infty}\alpha_k|\Theta(x^k)|<\infty~~\text{and}~~	\sum_{k=0}^{\infty}\alpha_k\|t^k\|^2<\infty.
		\end{equation}
	\end{lem}
	\begin{proof}
		By step 4 Algorithm 1, we know 	\begin{eqnarray*}
			\Phi_j(x^{k+1})\leq \Phi_j(x^k)+\alpha_k \beta \Phi^*_j(x^k,t^k),~~\forall~j\in \Lambda, \label{amj1}
		\end{eqnarray*}
		where  $\Phi^*_j(x^k,t^k) = \displaystyle \max_{i \in \bar \Lambda} \{f_j(x^k,\xi_i) + \nabla f_j(x^k,\xi_i)^Tt^k\}- \Phi_j(a^k)$.
	\end{proof}
	By taking the sum $k=0$ to $N$ we get
	\begin{align*}
		\Phi_j(x^{N+1})&\leq \Phi_j(x^0)+ \beta \sum_{k=0}^{N}\alpha_k\Phi^*_j(x^k,t^k)\\
		&\leq \Phi_j(x^0)+ \beta \sum_{k=0}^{N}\alpha_k\displaystyle \max_{i \in \bar \Lambda} \{f_j(x^k,\xi_i) + \nabla f_j(x^k,\xi_i)^Tt^k\}- \Phi_j(x^k)\\
		&\leq \Phi_j(x^0)+ \beta \sum_{k=0}^{N}\alpha_k\bigg\{\displaystyle \max_{j\in \Lambda}\max_{i \in \bar \Lambda} \{f_j(x^k,\xi_i) + \nabla f_j(x^k,\xi_i)^Tt^k- \Phi_j(a^k)\}+\frac{\|t^k\|^2}{2}-\frac{\|t^k\|^2}{2}\bigg\}\\
		&\leq \Phi_j(x^0)+ \beta \sum_{k=0}^{N}\alpha_k\bigg\{ \Theta(x^k)-\frac{\|t^k\|^2}{2}\bigg\}.
	\end{align*}
	The above inequality can be written as 
	\begin{equation*}
		\Phi_j(x^{N+1})-\Phi_j(x^0) \leq \beta \sum_{k=0}^{N}\alpha_k\bigg\{ \Theta(x^k)-\frac{\|t^k\|^2}{2}\bigg\},
	\end{equation*}
	as $\hat{y_j}\leq \Phi_j(x^k)$ and $\Theta(x^k)<0$ for all $j$ and $k$ therefore 
	\begin{equation*}
		\frac{1}{\beta}(\hat{y_j}-\Phi_j(x^0))\leq \sum_{k=0}^{N}\alpha_k\bigg( \Theta(x^k)-\frac{\|t^k\|^2}{2}\bigg),
	\end{equation*}
	also  \begin{equation*}
		\frac{1}{\beta}(\Phi_j(x^0)-\hat{y_j})\geq \sum_{k=0}^{N}\alpha_k\bigg( -\Theta(x^k)+\frac{\|t^k\|^2}{2}\bigg),
	\end{equation*}
	which implies \begin{equation*}
		\sum_{k=0}^{N}\alpha_k\bigg( |\Theta(x^k)|+\frac{\|t^k\|^2}{2}\bigg)\leq\frac{1}{\beta}(\Phi_j(x^0)-\hat{y_j}),~~\because~\theta(x^k)<0.
	\end{equation*}
	Since the right-hand side of the above inequality is finite and inequality holds for any positive integer $N$, then we get  
	\begin{equation*}
		\sum_{k=0}^{\infty}\alpha_k\bigg( |\Theta(x^k)|+\frac{\|t^k\|^2}{2}\bigg)<\infty.
	\end{equation*}
	The above inequality implies the result of this lemma.
\end{proof}
\begin{thm}
	Suppose that $\Phi$ is convex in component-wise sense (i.e., $\Phi$ is $\mathbb{R}^m-$ convex) and the Assumption 1 holds. Then any sequence produced by Algorithm \ref{algo1} converges to a WPOS $x^*\in \mathbb{R}^n.$ 
\end{thm}
\begin{proof}
	Since by Algorithm \ref{algo1}, $\{\Phi(x^k)\}$ is a component-wise decreasing sequence, then by assumption, there exists $\tilde{x}\in \mathbb{R}^n$ such that \begin{equation}\label{4.18}
		\Phi(\tilde{x})\leq \Phi(x^k)~\text{for all}~k=0,1,2\ldots.
	\end{equation}
	It is observed that $0<\alpha_k\leq1$ for all $k,$ so 
	\begin{align*}
		\|x^{k+1}-x^k\|^2&\leq\frac{1}{\alpha_k}\|x^{k+1}-x^k\|^2~~\text{for all} ~k=0,1,2,\ldots\\
		&\leq\frac{1}{\alpha_k}\|\alpha_k t^k\|^2=\alpha_k\| t^k\|^2~~\text{for all}~k=0,1,2,\ldots~ (\because x^{k+1}=x^k+\alpha_kt^k).
	\end{align*}
	Therefore, by above inequality, (\ref{4.18}), and Lemma (\ref{4.5}) we obtained
	$$\sum_{k=0}^{\infty}\|x^{k+1}-x^k\|^2\leq \sum_{k=0}^{\infty} \alpha_k\| t^k\|^2<\infty.$$ Thus, \begin{equation}\label{4.19}
		\sum_{k=0}^{\infty}\|x^{k+1}-x^k\|^2<\infty.
	\end{equation}
	Let us define $\tilde{L}=\{x\in \mathbb{R}^n: \Phi(x)\leq \Phi(x^k),~k=0,1,2,\ldots\}.$ By the component-wise convexity of $\Phi$ and Lemma \ref{4.3}, for any $x\in\tilde{L}$ we have \begin{equation*}
		\|x-x^{k+1}\|^2\leq \|x-x^{k}\|^2+ 	\|x^k-x^{k+1}\|^2~~\text{for all}~ k=0,1,2
		\ldots.
	\end{equation*}
	As $\tilde{L}$ is non empty because $\tilde{x}\in\tilde{L}$, by (\ref{4.19}) and the above inequality, it follows that $\{x^k\}$ is quasi-Fejer convergent to the set $\tilde{L}.$ Then by Theorem \ref{qfc}, $\{x^k\}$ is bounded and hence $\{x^k\}$ has an accumulation point. Let $x^*$ be one of them. Then by Lemma \ref{4.4}, $x^*\in\tilde{L}.$ Then by Theorem \ref{qfc}, we observe that $\{x^k\}$ converges to $x^*$. Therefore, by Theorem \ref{3.4}, $x^*$ is a critical point, and hence $\mathbb{R}^m-$ convexity implies that $x^*$ is a weak Pareto optimal solution for $\Phi.$
\end{proof}
\section{\textbf{Convergence rate of steepest descent method for UMOP}}\label{sec5}
Before establishing the convergence rate, we first present some results related to subproblem \eqref{eq31} which can be written as a $a=a^k$
\begin{eqnarray*}
	P(x^k):  \min_{t,\rho}~~~~\rho+\frac{1}{2}\|t\|^2&\\
	\text{s. t.}~~ f_j(x^k,\xi_i)+\nabla f_j(x^k,\xi_i)^Tt -\Phi_j(x^k)\leq &\rho ~\forall ~i\in \bar{\Lambda}~~j \in \Lambda.
\end{eqnarray*}
If we assume $\lambda^k_{ij}$ be the KKT optimal multiplier of the problem then the KKT optimality condition is given as follows:
\begin{eqnarray}
	\sum\limits_{j\in \Lambda}  \sum\limits_{i\in \bar{\Lambda}} \lambda^k_{ij} &=&1,\label{kkt1}\\
	t +\sum\limits_{j\in \Lambda} \sum\limits_{i\in \bar{\Lambda}} \lambda^k_{ij} \nabla f_j(x^k,\xi_i)&=&0,\label{kkt2}\\
	f_j(x^k,\xi_i)+\nabla f_j(x^k,\xi_i)^Tv -\Phi_j(x^k) -\rho&\leq& 0,~ \forall~ i\in \bar{\Lambda}~\forall~~j\in\Lambda, \label{kkt3}\\
	\lambda^k_{ij}\geq 0,~~~\lambda^k_{ij}\left(f_j(x^k,\xi_i)+\nabla f_j(x^k,\xi_i)^Tv  -\Phi_j(x^k)-\rho\right)&=&0, ~~\forall~~~ i\in \bar{\Lambda}~ ~\forall~ j\in\Lambda.\label{kkt4}
\end{eqnarray}
By (\ref{kkt2}),
\begin{equation}\label{4***}
	t(x^k) = -	\sum\limits_{j\in \Lambda}  \sum\limits_{i\in \bar \Lambda } \lambda^k_{ij} \nabla f_j(x^k,\xi_i).
\end{equation}
From (\ref{kkt4}), we have 
\begin{align}
	\rho&=\sum\limits_{j\in \Lambda}  \sum\limits_{i\in \bar{\Lambda}} \lambda^k_{ij}\big(f_j(x^k,\xi_i)+\nabla f_j(x^k,\xi_i)^Tt^k  -\Phi_j(x^k)\big)~~~~\because~t(x^k)=t^k\nonumber\\
	&\leq\sum\limits_{j\in \Lambda}  \sum\limits_{i\in \bar{\Lambda}} \lambda^k_{ij}\nabla f_j(x^k,\xi_i)^Tt^k \nonumber\\
	&=-\|t^k\|^2~~~\text{by}~~(\ref{4***}) \label{5.6}.
\end{align}
Now we make some assumptions to establish the convergence rate.
\begin{enumerate}[({As.}1)]
	\item\label{as1} $f_j(x,\xi_i)$ continuously differentiable and convex for each $j$ and $i.$
	\item\label{as2}For each $j$ and $i,$ $f_j(x,\xi_i)$ is $\gamma_i$-smooth, i.e., $\|\nabla f_j(x,\xi_i)-\nabla f_j(y,\xi_i)\|\leq \gamma_i\|x-y\|$ for every $x,y \in \mathbb{R}^n$ and $i\in\bar{\Lambda}.$
	\item\label{as3} At a given point $x^0,$ the objective function $\Phi$ has a bounded level set $L^S=\{x:\Phi(x)\leq \Phi(x^0)\}.$
\end{enumerate}
\begin{lem}\label{lem5.1} Let $\{x^k\}$ be a sequence generated by Algorithm \ref{algo1} with constant step size $\alpha_{k}=\alpha.$ If we assume As.\ref{as1}-As.\ref{as3} hold and $0<\alpha<\frac{2}{\gamma}.$ Then $\{\Phi_j(x^k)\}$ is decreasing i.e., $$\Phi_j(x^{k+1})\leq \Phi_j(x^k)-\alpha(1-\frac{\gamma}{2})\|t^k\|^2.$$ Moreover, for $k=0,1,2,\ldots,$ $$\sum_{k=0}^{\infty}\alpha(1-\frac{\gamma}{2}\alpha)\|{t^k}^2\|\leq \Phi_j(x^0)-y_j.$$
\end{lem} 
\begin{proof}
	Since $f_j(x,\xi_i)$ is $\gamma_i$ smooth then by iterative scheme of Algorithm\ref{algo1}, we have
	\begin{align*}
		f_j(x^{k+1},\xi_i)&\leq f_j(x^k,\xi_i)+\nabla f_j(x^k,\xi_i)^T(x^{k+1}-x^k)+\frac{\gamma}{2}\|x^{k+1}-x^k\|^2\\
		&\leq f_j(x^k,\xi_i)+\alpha\max_{j\in \Lambda}\max_{i \in \bar \Lambda}\{f_j(x^k,\xi_i)+\nabla f_j(x^k,\xi_i)^T t^k  -\Phi_j(x^k)\}+\alpha^2\frac{\gamma}{2}\|t^k\|^2\\
		&\leq f_j(x^k,\xi_i)-\alpha\|t^k\|^2+\alpha^2\frac{\gamma}{2}\|t^k\|^2,~~\text{by definition of}~\rho~\text{and since}~\rho\leq-\|t^k\|^2\\
		&\leq f_j(x^k,\xi_i)-\alpha(1-\frac{\gamma}{2}\alpha)\|t^k\|^2, ~\text{for each }~i~\text{and}~j\\
		\max_{i \in \bar \Lambda}	f_j(x^{k+1},\xi_i)&	\leq\max_{i \in \bar \Lambda}f_j(x^k,\xi_i)-\alpha(1-\frac{\gamma}{2}\alpha)\|t^k\|^2, ~\text{for each }~j.
	\end{align*}
	Therefore,
	\begin{align}\label{5.7}
		\Phi_j(x^{k+1})\leq \Phi_j(x^k)-\alpha(1-\frac{\gamma}{2}\alpha)\|t^k\|^2, ~\text{for each }~j.
	\end{align} 
	Since $0<\alpha<\frac{2}{\gamma}$, then last inequality (\ref{5.7}) implies that $\{\Phi_j(x^k)\}$ is decreasing sequence for each $j$ and hence each term of $\{x^k\}$ lies in the bounded level set $L^S.$ Due to the continuity of $\Phi_j,$ $\{\Phi_j(x^k)\}$ is bounded for each $j.$ Therefore, $\{\Phi_j(x^k)\}$ is convergent for each $j.$ Now, let us assume $\{\Phi_j(x^k)\}$ converges to $y_j.$ Now taking the summation $k=0$ to $\infty$ in inequality (\ref{5.7}), we get 
	\begin{align}\label{5.8}
		\sum_{k=0}^{\infty}\alpha(1-\frac{\gamma}{2}\alpha)\|t^k\|^2\leq \Phi_j(x^0)-y_j< \infty.
	\end{align}
	Hence, (\ref{5.7}) and (\ref{5.8}) show the proof of the Lemma.
\end{proof}
From (\ref{5.7}), it follows that for any $\delta>0,$ there exists a $k^0$ (positive integer) such that \begin{equation}\label{5.9}
	\sum_{k=k^0}^{\infty}\alpha^2(\alpha\gamma -1)\|t^k\|^2<\frac{1}{2}\delta, 
\end{equation}
in addition there also exists $k^N$ for this $\delta$ such that $\|x^{k_r}-a^*\|\leq\frac{\gamma}{2}$ for any $k_r\geq k^N,$ where $\{x^{k_r}\}$ is a subsequence of $\{x^k\}$ which converges to $x^*.$
This discussion will help to prove the next lemma.
\begin{lem}\label{lem5.2}
	Let $\{x^k\}$ be a sequence generated by Algorithm \ref{algo1}. If assumption \ref{as1}-\ref{as3} hold. Then 
	\begin{align}\label{5.10}
		\|x^{k+1}-x^*\|^2\leq\|x^k-x^*|+\alpha^2(\alpha\gamma -1)\|t^k\|^2. 
	\end{align}
\end{lem}
\begin{proof}
	In view of the discussion of solution of $P(a^k)$ and Algorithm \ref{algo1}, we obtained 
	\begin{align*}
		\|x^{k+1}-x^*\|^2&=\|x^k+\alpha t^k-x^*\|^2\\
		&\leq\|x^{k}-x^*\|^2-2\alpha\sum\limits_{j\in \Lambda}  \sum\limits_{i\in \bar \Lambda } \lambda^k_{ij} \nabla f_j(x^k,\xi_i)^T(x^k-x^*)+\alpha^2\|t^k\|^2.
	\end{align*}
	Due to the fact that $f_j(x^k,\xi_i)$ is convex for each $j$ and $i.$ We get \begin{align*}
		\|x^{k+1}-x^*\|^2
		&\leq\|x^{k}-x^*\|^2+2\alpha\bigg(\sum\limits_{j\in \Lambda}  \sum\limits_{i\in \bar \Lambda} \lambda^k_{ij} f_j(x^*,\xi_i)^T -\sum\limits_{j\in \Lambda}  \sum\limits_{i\in \bar \Lambda } \lambda^k_{ij}f_j(x^k,\xi_i)\bigg)+\alpha^2\|t^k\|^2.
	\end{align*}
	As we know by Lemma \ref{lem5.1}, $\{f_j(x^k)\}$ and $\{\Phi_j(x^k)\}$ both are decreasing sequence, then we obtained 
	\begin{align*}
		\|x^{k+1}-x^*\|^2
		&\leq\|x^{k}-x^*\|^2+2\alpha\bigg(\sum\limits_{j\in \Lambda}  \sum\limits_{i\in \bar \Lambda} \lambda^k_{ij} \nabla f_j( x^{k+1},\xi_i)-\sum\limits_{j\in \Lambda}  \sum\limits_{i\in \bar \Lambda } \lambda^k_{ij} f_j(x^k,\xi_i) \bigg)+\alpha^2\|t^k\|^2.
	\end{align*}
	By the smoothness of $\gamma$ (where $\gamma=\displaystyle\max_{j\in \Lambda}\gamma_j$) $\sum\limits_{i\in \bar \Lambda} \lambda^k_{ij} \nabla f_j( x,\xi_i)$ and decent property of $\Phi_j$ as well as $f_{j}$ we get
	\begin{align*}
		\|x^{k+1}-x^*\|^2
		&\leq\|x^{k}-x^*\|^2+2\alpha\bigg(\sum\limits_{j\in \Lambda}  \sum\limits_{i\in \bar \Lambda} \lambda^k_{ij} \nabla f_j( x^{k},\xi_i)^T(x^{k+1}-x^k)+\frac{\gamma}{2}\|x^{k+1}-x^k\|^2 \bigg)+\alpha^2\|t^k\|^2\\
		&=\|x^{k}-x^*\|^2+2\alpha\bigg(-\alpha\|t^k\|^2+\frac{\gamma}{2}\alpha^2\|t^k\|^2 \bigg)+\alpha^2\|t^k\|^2\\
		&=\|x^{k}-x^*\|^2-2\alpha^2\|t^k\|^2+\gamma\alpha^3\|t^k\|^2 +\alpha^2\|t^k\|^2\\
		&=\|x^{k}-x^*\|^2+(\gamma\alpha-1)\|t^k\|^2 \alpha^2
	\end{align*}
	the last inequality is the same as the inequality (\ref{5.10}).\\
\end{proof}
In the following theorem, we present that the sequence generated by Algorithm \ref{algo1} converges to the solution with a linear rate.
\begin{thm}
	Let $f_{j}(x,\xi_i)$ be continuously differentiable and strongly convex function for each $j\in\Lambda$ and $i\in\bar{\Lambda}.$ Let $\{x^k\}$ be any sequence generated by Algorithm \ref{algo1} with constant step length size $\alpha=\alpha_k$ for each $k$ such that  $0<\alpha\leq\frac{1}{\gamma}$. If As.\ref{as2} and As.\ref{as3} hold then $\{x^k\}$ converges to $x^*$ linearly more precisely $$\|x^{k+1}-x^k\|^2\leq(1-\beta\alpha)\|x^k-x^*\|^2$$
\end{thm}
\begin{proof}
	Set $r=\frac{1}{\alpha}.$ If $x^k$ is not a WPOS, the process of finding $x^{k+1}$ on the basis $\{x^k\}$ by using Algorithm \ref{algo1} Step 2, Step 3, and Step 4 with constant step size is given by 
	\begin{eqnarray}
		x^{k+1}&=&\underset{x\in R^n}{\mathrm{argmin}}\displaystyle \max_{j\in\Lambda}\max_{i\in \bar{\Lambda}}\{f_j(x,\xi_i)+\nabla f_j(x,\xi_i)^T(x-x^k)-\Phi_j(x)\}+\frac{r}{2}\|x-x^k\|^2.
	\end{eqnarray} 
	Now, we define a function 
	$$G^k(x)=\displaystyle\max_{j\in \Lambda}\max_{i \in \bar \Lambda}\{f_j(x,\xi_i)-f_{j}(x^k,\xi_i)\}.$$
	$G^k(x^k;x)$ is defined as the linearization of $G^k(x)$ at $x^k,$ which is given as follows: 
	$$G^k(x^k;x)=\displaystyle\max_{j\in \Lambda}\max_{i \in \bar \Lambda}\{f_j(x^k,\xi_i)+\nabla f_j(x^k,\xi_i)^T(x-x^k)-\Phi_j(x^k)\}.$$
	Now taking 
	$$G^k_r(x^k;x)=\displaystyle\max_{j\in \Lambda}\max_{i \in \bar \Lambda}\{f_j(x^k,\xi_i)+\nabla f_j(x^k,\xi_i)^T(x-x^k)-\Phi_j(x^k)\}+\frac{r}{2}\|x-x^k\|^2.$$
	Also, 
 $$x_{G^k(x^k;r)}=\underset{x\in \mathbb{R}^n}{\mathrm{argmin}}~G^k_r(x^k;x)=x^{k+1}.$$
	By $\alpha\leq\frac{1}{\gamma}$ and $r=\frac{1}{\alpha}$ it is implied that $\gamma\leq r.$ Then by (Proposition (we have to define)) we get 
	$$G^k(x)\geq G^k(x^{k+1})+r(x-x^k)^T(x^k-x^{k+1})+\frac{r}{2}\|x^k-x^{k+1}\|^2+\frac{\beta}{2}\|x^k-x\|^2.$$
	Now taking $x=x^*,$ $\alpha t^k=x^{k+1}-x^k$ and $r=\frac{1}{\alpha}$ in the above inequality, we get
	\begin{equation}\label{5.13}
		G^k(x^*)\geq G^k(x^{k+1})+(x^k-x^*)^Tt^k+\frac{\alpha}{2}\|t^k\|^2+\frac{\beta}{2}\|x^k-x^*\|^2.
	\end{equation}
Since by Lemma \ref{lem5.1} we have $f_{j}(x^*,\xi_i)\leq f_{j}(x^{k+1},\xi_i)$ for each $i,$ $j$ and $\Phi_j(x^*)\leq \Phi_j(x^{k+1}).$ Which implies $G^k(x^*)\leq G^k(x^{k+1})$, now we combine this relation with the above inequality (\ref{5.13}) then we obtained
\begin{equation*}\label{5.13}
		0\geq G^k(x^*)-G^k(x^{k+1})\geq(x^k-x^*)^Tt^k+\frac{\alpha}{2}\|t^k\|^2+\frac{\beta}{2}\|x^k-x^*\|^2,
	\end{equation*}
	therefore \begin{equation}\label{5.14}
		(x^k-x^*)^Tt^k\leq-\frac{\alpha}{2}\|t^k\|^2-\frac{\beta}{2}\|x^k-x^*\|^2.
	\end{equation}
	Also \begin{align}
		\|x^{k+1}-x^k\|^2&=\|x^k+\alpha t^k -x^* \|^2\nonumber\\
		&=\|x^k -x^* \|^2+2\alpha(x^k-x^*)^Tt^k+\alpha^2\| t^k\|^2\label{5.15}.
	\end{align}
	By inequality (\ref{5.14}) and (\ref{5.15}), we have 
	\begin{align*}\label{5.16}
		\|x^{k+1}-x^k\|^2&\leq\|x^k -x^* \|^2+2\alpha(-\frac{\alpha}{2}\|t^k\|^2-\frac{\beta}{2}\|x^k-x^*\|^2)+\alpha^2\| t^k\|^2\\
		&\leq\|x^k -x^* \|^2-\alpha^2\|t^k\|^2-\alpha\beta\|x^k-x^*\|^2+\alpha^2\| t^k\|^2\\
		&\leq\|x^k -x^* \|^2(1-\alpha\beta),
	\end{align*}
	and therefore $\|x^{k+1}-x^k\|^2\leq\|x^k -x^* \|^2(1-\alpha\beta),$ which is required inequality and hence $\{x^k\}$ is converges to $x^*$ with linear rate.
\end{proof}

\section{Conclusion}\label{sec6}
In this study, we have significantly advanced the understanding of the steepest descent method for uncertain multiobjective optimization problems by addressing critical gaps in our previous work \cite{kumar2024steepest}. While our earlier research focused on the local convergence analysis of the proposed steepest descent method for uncertain multiobjective optimization problems, it did not establish global convergence or analyze the rate of convergence.

We have successfully demonstrated global convergence and provided a detailed proof of linear rate of convergence analysis for the steepest descent method. These contributions offer a more robust theoretical foundation, bridging the gap between local and global convergence properties. By presenting comprehensive proofs and a nuanced analysis, our research not only clarifies the method's performance in various scenarios but also enhances its practical applicability.

Our findings reveal that the steepest descent method can be effectively utilized across a broader range of problem settings, offering improved efficiency and robustness. This advancement is expected to have a significant impact on the field of multiobjective optimization, facilitating more effective problem-solving strategies and practical implementations.

Future work could explore the application of these results to other optimization methods and investigate further improvements to the convergence properties. Overall, our study provides a deeper insight into the steepest descent method, contributing to its theoretical and practical advancement in optimization problems.
\section*{Acknowledgment:}
This research is supported by Govt. of India CSIR fellowship, Program No. 09/1174(0006)/2019-EMR-I, New Delhi.
\section*{Conflict of interest:}
The authors declare no conflicts of interest.


\bibliographystyle{sn-mathphys}
\bibliography{fullsdm_reference}

\providecommand{\noopsort}[1]{}\providecommand{\singleletter}[1]{#1}%

\begin{thebibliography}{45}
\ifx \bisbn   \undefined \def \bisbn  #1{ISBN #1}\fi
\ifx \binits  \undefined \def \binits#1{#1}\fi
\ifx \bauthor  \undefined \def \bauthor#1{#1}\fi
\ifx \batitle  \undefined \def \batitle#1{#1}\fi
\ifx \bjtitle  \undefined \def \bjtitle#1{#1}\fi
\ifx \bvolume  \undefined \def \bvolume#1{\textbf{#1}}\fi
\ifx \byear  \undefined \def \byear#1{#1}\fi
\ifx \bissue  \undefined \def \bissue#1{#1}\fi
\ifx \bfpage  \undefined \def \bfpage#1{#1}\fi
\ifx \blpage  \undefined \def \blpage #1{#1}\fi
\ifx \burl  \undefined \def \burl#1{\textsf{#1}}\fi
\ifx \doiurl  \undefined \def \doiurl#1{\url{https://doi.org/#1}}\fi
\ifx \betal  \undefined \def \betal{\textit{et al.}}\fi
\ifx \binstitute  \undefined \def \binstitute#1{#1}\fi
\ifx \binstitutionaled  \undefined \def \binstitutionaled#1{#1}\fi
\ifx \bctitle  \undefined \def \bctitle#1{#1}\fi
\ifx \beditor  \undefined \def \beditor#1{#1}\fi
\ifx \bpublisher  \undefined \def \bpublisher#1{#1}\fi
\ifx \bbtitle  \undefined \def \bbtitle#1{#1}\fi
\ifx \bedition  \undefined \def \bedition#1{#1}\fi
\ifx \bseriesno  \undefined \def \bseriesno#1{#1}\fi
\ifx \blocation  \undefined \def \blocation#1{#1}\fi
\ifx \bsertitle  \undefined \def \bsertitle#1{#1}\fi
\ifx \bsnm \undefined \def \bsnm#1{#1}\fi
\ifx \bsuffix \undefined \def \bsuffix#1{#1}\fi
\ifx \bparticle \undefined \def \bparticle#1{#1}\fi
\ifx \barticle \undefined \def \barticle#1{#1}\fi
\ifx \bconfdate \undefined \def \bconfdate #1{#1}\fi
\ifx \botherref \undefined \def \botherref #1{#1}\fi
\ifx \url \undefined \def \url#1{\textsf{#1}}\fi
\ifx \bchapter \undefined \def \bchapter#1{#1}\fi
\ifx \bbook \undefined \def \bbook#1{#1}\fi
\ifx \bcomment \undefined \def \bcomment#1{#1}\fi
\ifx \oauthor \undefined \def \oauthor#1{#1}\fi
\ifx \citeauthoryear \undefined \def \citeauthoryear#1{#1}\fi
\ifx \endbibitem  \undefined \def \endbibitem {}\fi
\ifx \bconflocation  \undefined \def \bconflocation#1{#1}\fi
\ifx \arxivurl  \undefined \def \arxivurl#1{\textsf{#1}}\fi
\csname PreBibitemsHook\endcsname

\bibitem{bhaskar2000applications}
\begin{barticle}
\bauthor{\bsnm{Bhaskar}, \binits{V.}},
\bauthor{\bsnm{Gupta}, \binits{S.K.}},
\bauthor{\bsnm{Ray}, \binits{A.K.}}:
\batitle{Applications of multiobjective optimization in chemical engineering}.
\bjtitle{Reviews in chemical engineering}
\bvolume{16}(\bissue{1}),
\bfpage{1}--\blpage{54}
(\byear{2000})
\end{barticle}
\endbibitem

\bibitem{stewart2008real}
\begin{botherref}
\oauthor{\bsnm{Stewart}, \binits{T.}},
\oauthor{\bsnm{Bandte}, \binits{O.}},
\oauthor{\bsnm{Braun}, \binits{H.}},
\oauthor{\bsnm{Chakraborti}, \binits{N.}},
\oauthor{\bsnm{Ehrgott}, \binits{M.}},
\oauthor{\bsnm{G{\"o}belt}, \binits{M.}},
\oauthor{\bsnm{Jin}, \binits{Y.}},
\oauthor{\bsnm{Nakayama}, \binits{H.}},
\oauthor{\bsnm{Poles}, \binits{S.}},
\oauthor{\bsnm{Di~Stefano}, \binits{D.}}:
Real-world applications of multiobjective optimization.
Multiobjective optimization: interactive and evolutionary approaches,
285--327
(2008)
\end{botherref}
\endbibitem

\bibitem{miettinen1999nonlinear}
\begin{botherref}
\oauthor{\bsnm{Miettinen}, \binits{K.}}:
Nonlinear multiobjective optimization.
Springer
(1999)
\end{botherref}
\endbibitem

\bibitem{deb2013multi}
\begin{bchapter}
\bauthor{\bsnm{Deb}, \binits{K.}},
\bauthor{\bsnm{Deb}, \binits{K.}}:
\bctitle{Multi-objective optimization}.
In: \bbtitle{Search Methodologies: Introductory Tutorials in Optimization and
  Decision Support Techniques},
pp. \bfpage{403}--\blpage{449}.
\bpublisher{Springer}, \blocation{???}
(\byear{2013})
\end{bchapter}
\endbibitem

\bibitem{ehrgott2005multicriteria}
\begin{botherref}
\oauthor{\bsnm{Ehrgott}, \binits{M.}}:
Multicriteria optimization.
Springer
(2005)
\end{botherref}
\endbibitem

\bibitem{fliege2000steepest}
\begin{barticle}
\bauthor{\bsnm{Fliege}, \binits{J.}},
\bauthor{\bsnm{Svaiter}, \binits{B.F.}}:
\batitle{Steepest descent methods for multicriteria optimization}.
\bjtitle{Mathematical methods of operations research}
\bvolume{51},
\bfpage{479}--\blpage{494}
(\byear{2000})
\end{barticle}
\endbibitem

\bibitem{fliege2009newton}
\begin{barticle}
\bauthor{\bsnm{Fliege}, \binits{J.}},
\bauthor{\bsnm{Drummond}, \binits{L.G.}},
\bauthor{\bsnm{Svaiter}, \binits{B.F.}}:
\batitle{Newton's method for multiobjective optimization}.
\bjtitle{SIAM Journal on Optimization}
\bvolume{20}(\bissue{2}),
\bfpage{602}--\blpage{626}
(\byear{2009})
\end{barticle}
\endbibitem

\bibitem{fukuda2014survey}
\begin{barticle}
\bauthor{\bsnm{Fukuda}, \binits{E.H.}},
\bauthor{\bsnm{Drummond}, \binits{L.M.G.}}:
\batitle{A survey on multiobjective descent methods}.
\bjtitle{Pesquisa Operacional}
\bvolume{34},
\bfpage{585}--\blpage{620}
(\byear{2014})
\end{barticle}
\endbibitem

\bibitem{ansary2015modified}
\begin{barticle}
\bauthor{\bsnm{Ansary}, \binits{M.A.}},
\bauthor{\bsnm{Panda}, \binits{G.}}:
\batitle{A modified quasi-newton method for vector optimization problem}.
\bjtitle{Optimization}
\bvolume{64}(\bissue{11}),
\bfpage{2289}--\blpage{2306}
(\byear{2015})
\end{barticle}
\endbibitem

\bibitem{lai2020q}
\begin{barticle}
\bauthor{\bsnm{Lai}, \binits{K.K.}},
\bauthor{\bsnm{Mishra}, \binits{S.K.}},
\bauthor{\bsnm{Ram}, \binits{B.}}:
\batitle{On q-quasi-newton's method for unconstrained multiobjective
  optimization problems}.
\bjtitle{Mathematics}
\bvolume{8}(\bissue{4}),
\bfpage{616}
(\byear{2020})
\end{barticle}
\endbibitem

\bibitem{birge2011introduction}
\begin{botherref}
\oauthor{\bsnm{Birge}, \binits{J.R.}},
\oauthor{\bsnm{Louveaux}, \binits{F.}}:
Introduction to stochastic programming.
Springer
(2011)
\end{botherref}
\endbibitem

\bibitem{doolittle2018robust}
\begin{barticle}
\bauthor{\bsnm{Doolittle}, \binits{E.K.}},
\bauthor{\bsnm{Kerivin}, \binits{H.L.}},
\bauthor{\bsnm{Wiecek}, \binits{M.M.}}:
\batitle{Robust multiobjective optimization with application to internet
  routing}.
\bjtitle{Annals of Operations Research}
\bvolume{271},
\bfpage{487}--\blpage{525}
(\byear{2018})
\end{barticle}
\endbibitem

\bibitem{chen2012including}
\begin{barticle}
\bauthor{\bsnm{Chen}, \binits{W.}},
\bauthor{\bsnm{Unkelbach}, \binits{J.}},
\bauthor{\bsnm{Trofimov}, \binits{A.}},
\bauthor{\bsnm{Madden}, \binits{T.}},
\bauthor{\bsnm{Kooy}, \binits{H.}},
\bauthor{\bsnm{Bortfeld}, \binits{T.}},
\bauthor{\bsnm{Craft}, \binits{D.}}:
\batitle{Including robustness in multi-criteria optimization for
  intensity-modulated proton therapy}.
\bjtitle{Physics in Medicine \& Biology}
\bvolume{57}(\bissue{3}),
\bfpage{591}
(\byear{2012})
\end{barticle}
\endbibitem

\bibitem{yu2013robust}
\begin{barticle}
\bauthor{\bsnm{Yu}, \binits{H.}},
\bauthor{\bsnm{Liu}, \binits{H.}}:
\batitle{Robust multiple objective game theory}.
\bjtitle{Journal of Optimization Theory and Applications}
\bvolume{159},
\bfpage{272}--\blpage{280}
(\byear{2013})
\end{barticle}
\endbibitem

\bibitem{fliege2014robust}
\begin{barticle}
\bauthor{\bsnm{Fliege}, \binits{J.}},
\bauthor{\bsnm{Werner}, \binits{R.}}:
\batitle{Robust multiobjective optimization \& applications in portfolio
  optimization}.
\bjtitle{European Journal of Operational Research}
\bvolume{234}(\bissue{2}),
\bfpage{422}--\blpage{433}
(\byear{2014})
\end{barticle}
\endbibitem

\bibitem{dellnitz2009computation}
\begin{barticle}
\bauthor{\bsnm{Dellnitz}, \binits{M.}},
\bauthor{\bsnm{Witting}, \binits{K.}}:
\batitle{Computation of robust pareto points}.
\bjtitle{International Journal of Computing Science and Mathematics}
\bvolume{2}(\bissue{3}),
\bfpage{243}--\blpage{266}
(\byear{2009})
\end{barticle}
\endbibitem

\bibitem{witting2013variational}
\begin{barticle}
\bauthor{\bsnm{Witting}, \binits{K.}},
\bauthor{\bsnm{Ober-Bl{\"o}baum}, \binits{S.}},
\bauthor{\bsnm{Dellnitz}, \binits{M.}}:
\batitle{A variational approach to define robustness for parametric
  multiobjective optimization problems}.
\bjtitle{Journal of Global Optimization}
\bvolume{57}(\bissue{2}),
\bfpage{331}--\blpage{345}
(\byear{2013})
\end{barticle}
\endbibitem

\bibitem{gunawan2005multi}
\begin{barticle}
\bauthor{\bsnm{Gunawan}, \binits{S.}},
\bauthor{\bsnm{Azarm}, \binits{S.}}:
\batitle{Multi-objective robust optimization using a sensitivity region
  concept}.
\bjtitle{Structural and Multidisciplinary Optimization}
\bvolume{29},
\bfpage{50}--\blpage{60}
(\byear{2005})
\end{barticle}
\endbibitem

\bibitem{barrico2006robustness}
\begin{bchapter}
\bauthor{\bsnm{Barrico}, \binits{C.}},
\bauthor{\bsnm{Antunes}, \binits{C.H.}}:
\bctitle{Robustness analysis in multi-objective optimization using a degree of
  robustness concept}.
In: \bbtitle{2006 IEEE International Conference on Evolutionary Computation},
pp. \bfpage{1887}--\blpage{1892}
(\byear{2006}).
\bcomment{IEEE}
\end{bchapter}
\endbibitem

\bibitem{goerigk2015algorithm}
\begin{botherref}
\oauthor{\bsnm{Goerigk}, \binits{M.}},
\oauthor{\bsnm{Sch{\"o}bel}, \binits{A.}}:
Algorithm engineering in robust optimization.
arXiv preprint arXiv:1505.04901
(2015)
\end{botherref}
\endbibitem

\bibitem{jeyakumar2012robust}
\begin{barticle}
\bauthor{\bsnm{Jeyakumar}, \binits{V.}},
\bauthor{\bsnm{Li}, \binits{G.}},
\bauthor{\bsnm{Lee}, \binits{G.}}:
\batitle{Robust duality for generalized convex programming problems under data
  uncertainty}.
\bjtitle{Nonlinear Analysis: Theory, Methods \& Applications}
\bvolume{75}(\bissue{3}),
\bfpage{1362}--\blpage{1373}
(\byear{2012})
\end{barticle}
\endbibitem

\bibitem{goberna2013robust}
\begin{barticle}
\bauthor{\bsnm{Goberna}, \binits{M.A.}},
\bauthor{\bsnm{Jeyakumar}, \binits{V.}},
\bauthor{\bsnm{Li}, \binits{G.}},
\bauthor{\bsnm{L{\'o}pez}, \binits{M.A.}}:
\batitle{Robust linear semi-infinite programming duality under uncertainty}.
\bjtitle{Mathematical Programming}
\bvolume{139}(\bissue{1}),
\bfpage{185}--\blpage{203}
(\byear{2013})
\end{barticle}
\endbibitem

\bibitem{ide2014relationship}
\begin{barticle}
\bauthor{\bsnm{Ide}, \binits{J.}},
\bauthor{\bsnm{K{\"o}bis}, \binits{E.}},
\bauthor{\bsnm{Kuroiwa}, \binits{D.}},
\bauthor{\bsnm{Sch{\"o}bel}, \binits{A.}},
\bauthor{\bsnm{Tammer}, \binits{C.}}:
\batitle{The relationship between multi-objective robustness concepts and
  set-valued optimization}.
\bjtitle{Fixed Point Theory and Applications}
\bvolume{2014},
\bfpage{1}--\blpage{20}
(\byear{2014})
\end{barticle}
\endbibitem

\bibitem{ehrgott2014minmax}
\begin{barticle}
\bauthor{\bsnm{Ehrgott}, \binits{M.}},
\bauthor{\bsnm{Ide}, \binits{J.}},
\bauthor{\bsnm{Sch{\"o}bel}, \binits{A.}}:
\batitle{Minmax robustness for multi-objective optimization problems}.
\bjtitle{European Journal of Operational Research}
\bvolume{239}(\bissue{1}),
\bfpage{17}--\blpage{31}
(\byear{2014})
\end{barticle}
\endbibitem

\bibitem{goberna2015robust}
\begin{barticle}
\bauthor{\bsnm{Goberna}, \binits{M.A.}},
\bauthor{\bsnm{Jeyakumar}, \binits{V.}},
\bauthor{\bsnm{Li}, \binits{G.}},
\bauthor{\bsnm{Vicente-P{\'e}rez}, \binits{J.}}:
\batitle{Robust solutions to multi-objective linear programs with uncertain
  data}.
\bjtitle{European Journal of Operational Research}
\bvolume{242}(\bissue{3}),
\bfpage{730}--\blpage{743}
(\byear{2015})
\end{barticle}
\endbibitem

\bibitem{chuong2016optimality}
\begin{barticle}
\bauthor{\bsnm{Chuong}, \binits{T.D.}}:
\batitle{Optimality and duality for robust multiobjective optimization
  problems}.
\bjtitle{Nonlinear Analysis}
\bvolume{134},
\bfpage{127}--\blpage{143}
(\byear{2016})
\end{barticle}
\endbibitem

\bibitem{lee2018optimality}
\begin{barticle}
\bauthor{\bsnm{Lee}, \binits{J.H.}},
\bauthor{\bsnm{Lee}, \binits{G.M.}}:
\batitle{On optimality conditions and duality theorems for robust semi-infinite
  multiobjective optimization problems}.
\bjtitle{Annals of Operations Research}
\bvolume{269}(\bissue{1}),
\bfpage{419}--\blpage{438}
(\byear{2018})
\end{barticle}
\endbibitem

\bibitem{kumar2024steepest}
\begin{botherref}
\oauthor{\bsnm{Kumar}, \binits{S.}},
\oauthor{\bsnm{Ansary}, \binits{M.A.T.}},
\oauthor{\bsnm{Mahato}, \binits{N.K.}},
\oauthor{\bsnm{Ghosh}, \binits{D.}}:
Steepest descent method for uncertain multiobjective optimization problems
  under finite uncertainty set.
Applicable Analysis,
1--22
(2024)
\end{botherref}
\endbibitem

\bibitem{shubham2023newton}
\begin{barticle}
\bauthor{\bsnm{Kumar}, \binits{S.}},
\bauthor{\bsnm{Ansary}, \binits{M.A.T.}},
\bauthor{\bsnm{Mahato}, \binits{N.K.}},
\bauthor{\bsnm{Ghosh}, \binits{D.}},
\bauthor{\bsnm{Shehu}, \binits{Y.}}:
\batitle{Newton's method for uncertain multiobjective optimization problems
  under finite uncertainty set}.
\bjtitle{Journal of Nonlinear and Variational Analysis}
\bvolume{7}(\bissue{5}),
\bfpage{785}--\blpage{809}
(\byear{2023})
\end{barticle}
\endbibitem

\bibitem{kumar2024modified}
\begin{botherref}
\oauthor{\bsnm{Kumar}, \binits{S.}},
\oauthor{\bsnm{Mahato}, \binits{N.K.}},
\oauthor{\bsnm{Ansary}, \binits{M.A.T.}},
\oauthor{\bsnm{Ghosh}, \binits{D.}},
\oauthor{\bsnm{Trean{\c{t}}{\u{a}}}, \binits{S.}}:
A modified quasi-newton method for uncertain multiobjective optimization
  problems under a finite uncertainty set.
Engineering Optimization,
1--30
(2024)
\end{botherref}
\endbibitem

\bibitem{drummond2005steepest}
\begin{barticle}
\bauthor{\bsnm{Drummond}, \binits{L.G.}},
\bauthor{\bsnm{Svaiter}, \binits{B.F.}}:
\batitle{A steepest descent method for vector optimization}.
\bjtitle{Journal of computational and applied mathematics}
\bvolume{175}(\bissue{2}),
\bfpage{395}--\blpage{414}
(\byear{2005})
\end{barticle}
\endbibitem

\bibitem{povalej2014quasi}
\begin{barticle}
\bauthor{\bsnm{Povalej}, \binits{{\v{Z}}.}}:
\batitle{Quasi-newton's method for multiobjective optimization}.
\bjtitle{Journal of Computational and Applied Mathematics}
\bvolume{255},
\bfpage{765}--\blpage{777}
(\byear{2014})
\end{barticle}
\endbibitem

\bibitem{mahdavi2020superlinearly}
\begin{barticle}
\bauthor{\bsnm{Mahdavi-Amiri}, \binits{N.}},
\bauthor{\bsnm{Salehi~Sadaghiani}, \binits{F.}}:
\batitle{A superlinearly convergent nonmonotone quasi-newton method for
  unconstrained multiobjective optimization}.
\bjtitle{Optimization Methods and Software}
\bvolume{35}(\bissue{6}),
\bfpage{1223}--\blpage{1247}
(\byear{2020})
\end{barticle}
\endbibitem

\bibitem{morovati2018quasi}
\begin{barticle}
\bauthor{\bsnm{Morovati}, \binits{V.}},
\bauthor{\bsnm{Basirzadeh}, \binits{H.}},
\bauthor{\bsnm{Pourkarimi}, \binits{L.}}:
\batitle{Quasi-newton methods for multiobjective optimization problems}.
\bjtitle{4OR}
\bvolume{16},
\bfpage{261}--\blpage{294}
(\byear{2018})
\end{barticle}
\endbibitem

\bibitem{qu2011quasi}
\begin{barticle}
\bauthor{\bsnm{Qu}, \binits{S.}},
\bauthor{\bsnm{Goh}, \binits{M.}},
\bauthor{\bsnm{Chan}, \binits{F.T.}}:
\batitle{Quasi-newton methods for solving multiobjective optimization}.
\bjtitle{Operations Research Letters}
\bvolume{39}(\bissue{5}),
\bfpage{397}--\blpage{399}
(\byear{2011})
\end{barticle}
\endbibitem

\bibitem{lucambio2018nonlinear}
\begin{barticle}
\bauthor{\bsnm{Lucambio~P{\'e}rez}, \binits{L.}},
\bauthor{\bsnm{Prudente}, \binits{L.}}:
\batitle{Nonlinear conjugate gradient methods for vector optimization}.
\bjtitle{SIAM Journal on Optimization}
\bvolume{28}(\bissue{3}),
\bfpage{2690}--\blpage{2720}
(\byear{2018})
\end{barticle}
\endbibitem

\bibitem{gonccalves2020extension}
\begin{barticle}
\bauthor{\bsnm{Gon{\c{c}}alves}, \binits{M.L.}},
\bauthor{\bsnm{Prudente}, \binits{L.}}:
\batitle{On the extension of the hager--zhang conjugate gradient method for
  vector optimization}.
\bjtitle{Computational Optimization and Applications}
\bvolume{76}(\bissue{3}),
\bfpage{889}--\blpage{916}
(\byear{2020})
\end{barticle}
\endbibitem

\bibitem{drummond2004projected}
\begin{barticle}
\bauthor{\bsnm{Drummond}, \binits{L.G.}},
\bauthor{\bsnm{Iusem}, \binits{A.N.}}:
\batitle{A projected gradient method for vector optimization problems}.
\bjtitle{Computational Optimization and applications}
\bvolume{28},
\bfpage{5}--\blpage{29}
(\byear{2004})
\end{barticle}
\endbibitem

\bibitem{fukuda2011convergence}
\begin{barticle}
\bauthor{\bsnm{Fukuda}, \binits{E.H.}},
\bauthor{\bsnm{Drummond}, \binits{L.G.}}:
\batitle{On the convergence of the projected gradient method for vector
  optimization}.
\bjtitle{Optimization}
\bvolume{60}(\bissue{8-9}),
\bfpage{1009}--\blpage{1021}
(\byear{2011})
\end{barticle}
\endbibitem

\bibitem{fukuda2013inexact}
\begin{barticle}
\bauthor{\bsnm{Fukuda}, \binits{E.H.}},
\bauthor{\bsnm{Gra{\~n}a~Drummond}, \binits{L.}}:
\batitle{Inexact projected gradient method for vector optimization}.
\bjtitle{Computational Optimization and Applications}
\bvolume{54},
\bfpage{473}--\blpage{493}
(\byear{2013})
\end{barticle}
\endbibitem

\bibitem{fazzio2019convergence}
\begin{barticle}
\bauthor{\bsnm{Fazzio}, \binits{N.S.}},
\bauthor{\bsnm{Schuverdt}, \binits{M.L.}}:
\batitle{Convergence analysis of a nonmonotone projected gradient method for
  multiobjective optimization problems}.
\bjtitle{Optimization Letters}
\bvolume{13},
\bfpage{1365}--\blpage{1379}
(\byear{2019})
\end{barticle}
\endbibitem

\bibitem{bonnel2005proximal}
\begin{barticle}
\bauthor{\bsnm{Bonnel}, \binits{H.}},
\bauthor{\bsnm{Iusem}, \binits{A.N.}},
\bauthor{\bsnm{Svaiter}, \binits{B.F.}}:
\batitle{Proximal methods in vector optimization}.
\bjtitle{SIAM Journal on Optimization}
\bvolume{15}(\bissue{4}),
\bfpage{953}--\blpage{970}
(\byear{2005})
\end{barticle}
\endbibitem

\bibitem{ceng2010hybrid}
\begin{barticle}
\bauthor{\bsnm{Ceng}, \binits{L.}},
\bauthor{\bsnm{Mordukhovich}, \binits{B.S.}},
\bauthor{\bsnm{Yao}, \binits{J.-C.}}:
\batitle{Hybrid approximate proximal method with auxiliary variational
  inequality for vector optimization}.
\bjtitle{Journal of optimization theory and applications}
\bvolume{146}(\bissue{2}),
\bfpage{267}--\blpage{303}
(\byear{2010})
\end{barticle}
\endbibitem

\bibitem{sun2006optimization}
\begin{botherref}
\oauthor{\bsnm{Sun}, \binits{W.}},
\oauthor{\bsnm{Yuan}, \binits{Y.-X.}}:
Optimization theory and methods: nonlinear programming.
Springer
(2006)
\end{botherref}
\endbibitem

\bibitem{shubham2024steepest}
\begin{botherref}
\oauthor{\bsnm{Kumar}, \binits{S.}},
\oauthor{\bsnm{Ansary}, \binits{M.A.T.}},
\oauthor{\bsnm{Mahato}, \binits{N.K.}},
\oauthor{\bsnm{Ghosh}, \binits{D.}}:
Steepest descent method for uncertain multiobjective optimization problems
  under finite uncertainty set.
Applicable Analysis,
1--25
(2024)
\end{botherref}
\endbibitem

\end{thebibliography}

\end{document}